\documentclass[11pt]{article}

\usepackage[margin=1in]{geometry}
\usepackage{amsmath,amssymb,amsthm,mathtools}
\usepackage{microtype}
\usepackage[hidelinks]{hyperref}

\newtheorem{theorem}{Theorem}
\newtheorem{proposition}{Proposition}
\newtheorem{remark}{Remark}

\newcommand{\R}{\mathbb{R}}
\newcommand{\B}{\mathrm{B}}
\newcommand{\dd}{\,d}

\title{H\"older Duality and Self-Dual Arc Length on \(L^p\) Unit Curves}
\author{Jacinto Joaqu\'in Arichabala Concha}
\date{}

\begin{document}
\maketitle

\begin{abstract}
We study arc length on the first-quadrant unit curve
\[
x^p+y^p=1,\qquad p>1,
\]
when length is measured with respect to an \(L^r\) norm. Denoting this
length by \(L(p,r)\), we prove the H\"older-duality identity
\[
L(p,r)=L(r^*,p^*),
\]
where \(p^*=p/(p-1)\) and \(r^*=r/(r-1)\).
The identity is interpreted as a special case of the symmetry of mixed
area (equivalently, of the anisotropic perimeter pairing).
For the distinguished dual metric \(r=p^*\), the length element agrees
with the areal angular element:
\[
L(p,p^*)=\int_\gamma (x\,dy-y\,dx)
=\frac1p\,\B\!\left(\frac1p,\frac1p\right).
\]
Thus the familiar coincidence between arc length and angular parameter
on the Euclidean unit circle is the self-dual case \(p=2\) of a more
general relation between \(L^p\) geometry and its dual metric.  As a
consequence, the Euclidean arc length of the \(L^p\) unit curve admits a
dual trigonometric representation.
\end{abstract}

\section{Introduction}

The Euclidean unit circle has the familiar property that, on the unit
circle, arc length and angular parameter coincide.  Equivalently, the
arc length of an oriented sector equals twice its swept area.  This
phenomenon is closely tied to the self-duality of the Euclidean norm.

We consider the first-quadrant arc of
\[
x^p+y^p=1,\qquad p>1,
\]
and measure its length using an \(L^r\) norm.  If \(L(p,r)\) denotes this
length, we show
\[
L(p,r)=L(r^*,p^*),
\]
where \(p^*\) and \(r^*\) are the H\"older conjugates.  In convex-geometric
language, this identity is a specialization of the symmetry of mixed
area.  We retain an elementary support-function proof because it makes
the interchange between the exponent defining the curve and the
exponent defining the metric explicit.

The dual case \(r=p^*\) is especially transparent: the \(L^{p^*}\) arc
length element on the \(L^p\) unit curve equals the differential
\(x\,dy-y\,dx\).  Hence dual arc length equals twice swept sectorial area.
For \(p=2\), geometry and metric are self-dual and the classical
Euclidean relation is recovered.

\section{Definition and an integral representation}

Let
\[
B_p=\{(x,y)\in\R^2: |x|^p+|y|^p\le 1\},
\qquad p>1,
\]
and let \(\gamma_p\) be the first-quadrant arc of \(\partial B_p\).
For \(r>1\), define
\[
L(p,r)=\int_{\gamma_p}\|d\gamma\|_r.
\]

A convenient parametrization is
\[
x(\theta)=\cos^{2/p}\theta,\qquad
y(\theta)=\sin^{2/p}\theta,\qquad
0\le\theta\le\frac{\pi}{2}.
\]
Differentiation gives
\[
\left|\frac{dx}{d\theta}\right|
=\frac{2}{p}\sin\theta\,\cos^{2/p-1}\theta,
\qquad
\left|\frac{dy}{d\theta}\right|
=\frac{2}{p}\cos\theta\,\sin^{2/p-1}\theta.
\]
Consequently,
\begin{align}
L(p,r)
&=\frac{2}{p}\int_0^{\pi/2}
\left(
\sin^r\theta\,\cos^{r(2/p-1)}\theta+
\cos^r\theta\,\sin^{r(2/p-1)}\theta
\right)^{1/r}\dd\theta \notag\\
&=\frac{2}{p}\int_0^{\pi/2}
(\sin\theta\cos\theta)^{2/p-1}
\left(
\sin^{2r/p^*}\theta+\cos^{2r/p^*}\theta
\right)^{1/r}\dd\theta.
\label{eq:integral}
\end{align}

\section{H\"older duality for arc length}

For \(q>1\), let \(h_q\) denote the support function of \(B_q\).  If
\(n(\phi)=(\cos\phi,\sin\phi)\), then
\[
h_q(\phi)=\|n(\phi)\|_{q^*}
=\left(\cos^{q^*}\phi+\sin^{q^*}\phi\right)^{1/q^*}.
\]

\begin{theorem}[H\"older duality]
For all \(p,r>1\),
\[
\boxed{L(p,r)=L(r^*,p^*).}
\]
\end{theorem}

\begin{proof}
Parametrize the first-quadrant arc of \(\partial B_p\) by its Euclidean
outer-normal angle \(\phi\).  On the open quadrant the Euclidean arc
length element is
\[
ds_p=(h_p+h_p'')\dd\phi.
\]
The positively oriented tangent is a rotation of the unit normal, hence
\[
\|d\gamma_p\|_r
=\|(-\sin\phi,\cos\phi)\|_r\,ds_p
=h_{r^*}(\phi)\,ds_p.
\]
Therefore
\begin{equation}
L(p,r)=\int_0^{\pi/2}
h_{r^*}(\phi)\bigl(h_p(\phi)+h_p''(\phi)\bigr)\dd\phi.
\label{eq:pairing}
\end{equation}

A minor endpoint issue deserves explicit mention.  If
\(s=q^*>1\), then
\[
h_q(\phi)=\bigl(\cos^s\phi+\sin^s\phi\bigr)^{1/s}
\]
satisfies
\[
h_q'(0)=h_q'(\pi/2)=0.
\]
Its second derivative may be unbounded at an endpoint when \(1<s<2\),
but its possible singularity is of order \(\phi^{s-2}\), hence is
integrable.  Thus the following integrations by parts may be performed
first on \([\varepsilon,\pi/2-\varepsilon]\) and then passed to the limit
as \(\varepsilon\downarrow0\).  The boundary term vanishes, and
\[
\int_0^{\pi/2}h_{r^*}h_p''\dd\phi
=
\int_0^{\pi/2}h_p h_{r^*}''\dd\phi.
\]
Using \eqref{eq:pairing},
\[
L(p,r)
=\int_0^{\pi/2}
h_p(\phi)\bigl(h_{r^*}(\phi)+h_{r^*}''(\phi)\bigr)\dd\phi.
\]
The last expression is the \(L^{p^*}\) length of the first-quadrant arc
of \(\partial B_{r^*}\), since the support function associated with the
dual metric \(L^{p^*}\) is \(h_p\).  Hence
\[
L(p,r)=L(r^*,p^*).
\]
\end{proof}

\begin{remark}[Mixed-area interpretation]
For planar convex bodies \(K,L\), the mixed area admits the support
function representation
\[
V(K,L)=\frac12\int_{\mathbb S^1}h_L\,dS_K,
\]
and is symmetric: \(V(K,L)=V(L,K)\).  By the fourfold symmetry of
\(L^p\) balls, the preceding theorem is precisely this symmetry applied
to \(K=B_p\) and \(L=B_{r^*}\), restricted to one quadrant.  Thus the
H\"older-duality identity is a concrete \(L^p\)-geometric manifestation
of the classical symmetry of mixed area.
\end{remark}

\section{The dual metric and the areal angular parameter}

Set \(r=p^*\) in \eqref{eq:integral}.  Since
\[
\sin^2\theta+\cos^2\theta=1,
\]
we obtain
\begin{equation}
L(p,p^*)
=\frac{2}{p}\int_0^{\pi/2}
(\sin\theta\cos\theta)^{2/p-1}\dd\theta.
\label{eq:selfdual}
\end{equation}

Define the areal angular parameter along \(\gamma_p\) by
\[
\alpha=\int_{\gamma_p}(x\,dy-y\,dx).
\]
It is twice the signed area swept by the radius vector.  Directly from
the parametrization,
\[
x\,dy-y\,dx
=\frac{2}{p}
(\sin\theta\cos\theta)^{2/p-1}\dd\theta.
\]
Comparison with \eqref{eq:selfdual} gives the following.

\begin{proposition}
For every \(p>1\),
\[
\boxed{L(p,p^*)=\alpha.}
\]
Equivalently, the \(L^{p^*}\) length of the first-quadrant arc of
\(\partial B_p\) equals twice the area of the corresponding quadrant of
\(B_p\).
\end{proposition}

Using
\[
\int_0^{\pi/2}\sin^{a-1}\theta\cos^{b-1}\theta\dd\theta
=\frac12\B\!\left(\frac a2,\frac b2\right),
\]
we obtain
\begin{equation}
\boxed{
L(p,p^*)
=\frac1p\B\!\left(\frac1p,\frac1p\right)
=\frac1p\frac{\Gamma(1/p)^2}{\Gamma(2/p)}.
}
\end{equation}

For \(p=2\), one has \(p^*=2\), so the curve and the metric are
self-dual.  The preceding identity reduces to the usual coincidence of
unit-circle arc length with angular parameter.  Over the entire unit
circle it also gives the familiar relation
\[
\operatorname{Perimeter}(B_2)=2\,\operatorname{Area}(B_2).
\]

\section{Euclidean arc length as a dual representation}

Taking \(r=2\) in the duality theorem and using \(2^*=2\),
\[
L(p,2)=L(2,p^*).
\]
Parametrizing the Euclidean unit circle by
\(x=\cos\theta,\ y=\sin\theta\), we obtain
\[
\boxed{
L(p,2)=
\int_0^{\pi/2}
\left(\cos^{p^*}\theta+\sin^{p^*}\theta\right)^{1/p^*}\dd\theta.
}
\]

\section{Concluding remarks}

The two-parameter length functional \(L(p,r)\) makes explicit an
interchange between the exponent defining an \(L^p\) unit curve and the
exponent defining the metric in which that curve is measured.  Its
H\"older-duality symmetry is naturally explained by mixed area.

The case \(r=p^*\) is distinguished: dual arc length becomes exactly an
areal angular parameter.  In this sense the Euclidean unit circle is the
self-dual member of a broader family in which the classical equality
between arc length and angular parameter persists after replacing the
ambient metric by the H\"older-dual metric.

A further question, left open here, concerns the diagonal quantity
\(p\mapsto L(p,p)\).  Numerical evidence suggests that it may attain its
global minimum at \(p=2\).  We do not use this observation as a theorem
in the present paper.

\end{document}